 \documentclass[final,3p,times]{elsarticle}

\usepackage{amssymb}
 \usepackage{amsthm}
 \usepackage[T1]{fontenc}
 
\usepackage[ansinew]{inputenc} 
\usepackage{mathrsfs}
\usepackage{euscript}
\usepackage{natbib}

\newtheorem{theorem}{Theorem}
\newtheorem{definition}{Definition}
\newtheorem{lemma}{Lemma}
\newtheorem{proposition}{Proposition}
\newdefinition{rmk}{ Remark}
\newproof{pf}{Proof}

\begin{document}

\begin{frontmatter}



\title{ Critical exponent for damped wave equations with nonlinear memory}

 \author[label1,label2]{Ahmad Z. FINO}
 \address[label1]{Laboratoire de math\'ematiques appliqu\'ees, UMR CNRS 5142, Universit\'e de Pau et des Pays
de l'Adour, 64000 Pau, France}
\address[label2]{LaMA-Liban,
Lebanese University, P.O. Box 37 Tripoli, Lebanon}


 \ead{ahmad.fino01@gmail.com}
 \begin{abstract}
We consider the Cauchy problem in $\mathbb{R}^n,$ $n\geq 1,$ for a semilinear damped wave equation with nonlinear memory. Global existence and asymptotic behavior as $t\rightarrow\infty$ of small data solutions have been established in the case when $1\leq n\leq3.$ Moreover, we derive a blow-up result under some positive data in any dimensional space.
\end{abstract}

\begin{keyword}
Nonlinear damped wave equation\sep Global existence \sep  Blow-up \sep Critical exponent \sep Large time asymptotic behavior

\MSC[2010] 35L15 \sep 35L70 \sep 35B33 \sep 34B44
\end{keyword}

\end{frontmatter}


\section{Introduction}
\setcounter{equation}{0}

This paper concerns with the Cauchy problem for the damped wave equation with nonlinear memory
\begin{equation}\label{eq1}
\left\{\begin{array}{ll}
\,\, \displaystyle {u_{tt}-\Delta
u + u_t=\int_0^t(t-s)^{-\gamma}|u(s)|^p\,ds} &\displaystyle {t>0,x\in {\mathbb{R}^n},}\\
{}\\
\displaystyle{u(0,x)= u_0(x),\;\;u_t(0,x)= u_1(x)\qquad\qquad}&\displaystyle{x\in {\mathbb{R}^n},}
\end{array}
\right.
\end{equation} 
where the unknown function $u$ is real-valued, $n\geq 1,$ $0<\gamma<1$ and $p>1.$ Throughout this paper, we assume that 
\begin{equation}\label{condition1}
(u_0,u_1)\in H^1(\mathbb{R}^n)\times L^2(\mathbb{R}^n)
\end{equation}
 and
\begin{equation}\label{condition2}
\mbox{supp}u_i\subset B(K):=\{x\in\mathbb{R}^n:\;|x|< K\},\quad K>0,\;i=0,1.
\end{equation}
For the simplicity of notations, $\|\cdotp\|_q$ and $\|\cdotp\|_{H^1}$ $(1\leq q\leq \infty)$ stand for the usual $L^q(\mathbb{R}^n)$-norm and $H^1(\mathbb{R}^n)$-norm, respectively.

The nonlinear nonlocal term can be considered as an approximation of the classical semilinear damped
wave equation
$$u_{tt}-\Delta
u+u_t=|u(t)|^p\, $$
since the limit
$$ \lim_{\gamma \rightarrow 1} \frac{1}{\Gamma(1- \gamma)} s_+^{-\gamma} = \delta(s)$$
exists in distribution sense, where $\Gamma$ is the Euler gamma function.

\noindent It is clear that this nonlinear term involves memory type selfinteraction and can be considered  as Riemann-Liouville integral operator
$$
J^\alpha_{a|t}g(t):=\frac{1}{\Gamma(\alpha)}\int_a^t(t-s)^{\alpha-1}g(s)\,ds
$$
introduced with $a=-\infty$
by Liouville in 1832 and with $a=0$
 by Riemann in 1876 (see Chapter V in \cite{DB}).
Therefore, (\ref{eq1}) takes the form

\begin{equation}\label{A.2++}
     u_{tt}-\Delta
u+u_t=J^{\alpha}_{0|t}\left(|u|^p\right)(t),
\end{equation}
where $\alpha = 1 - \gamma.$

In recent years, questions of global existence and blow-up of solutions for nonlinear hyperbolic equations with a damping term have been studied by many mathematicians, see \cite{Ikehatatanizawa,Karch,Matsumura,Nishiharazhao,Todorovayordanov} and the references therein. To focus on our motivation, we shall mention below only some results related to Todorova and Yordanov \cite{Todorovayordanov}. For the Cauchy problem for the semilinear damped wave equation with the forcing term
\begin{equation}\label{classicalcase}
u_{tt}-\Delta u + u_t=|u|^p,\qquad u(0)=u_0,\quad u_t(0)=u_1,
\end{equation}
it has been conjectured that the damped wave equation has the diffuse structure as $t\rightarrow\infty$ (see e.g. \cite{Belloutfriedman,Li}). This suggests that problem $(\ref{classicalcase})$ should have $p_c(n):=1+{2}/{n}$ as critical exponent which is called the Fujita exponent named after Fujita \cite{Fujita}, in general space dimension. Indeed, Todorova and Yordanov \cite{Todorovayordanov} have showed that the critical exponent is exactly $p_c(n),$ that is, if $p>p_c(n)$ then all small initial data solutions of $(\ref{classicalcase})$ are global, while if $1<p<p_c(n)$ then all solutions of $(\ref{classicalcase})$ with initial data having positive average value blow-up in finite time regardless of the smallness of the initial data. Moreover, they showed that in the case of $p>p_c(n),$ the support of the solution of $(\ref{classicalcase})$ is strongly suppressed by the damping, so that the solution is concentrated in a ball much smaller than $|x|<t+K,$ namely
$$\|Du(t,\cdotp)\|_{L^2(\mathbb{R}^{n}\setminus B(t^{{1/2}+\delta}))}=\mathcal{O}(e^{-t^{2{\delta}/{4}}}),\quad\mbox{as}\;\; t\rightarrow\infty,$$
where $D:=(\partial_t,\nabla_x).$ Furthermore, they proved that the total energy of the solutions of $(\ref{classicalcase})$ decays at the rate of the linear equation, namely
$$\|Du(t,\cdotp)\|_{L^2(\mathbb{R}^n)}=\mathcal{O}(t^{-n/4-1/2}),\quad \mbox{as}\;\;t\rightarrow\infty.$$

Our goal is to apply the above properties founded by Todorova and Yordanov to our problem $(\ref{eq1})$ with the same assumptions on the initial data. The method used to prove the global existence is inspired from the weighted energy method developed in \cite{Todorovayordanov}. On the other hand, the test function method (see \cite{Finogeorgiev,FinoKarch,Finokirane, GK, KLT, PM1,PM2, Zhang} and the references therein) is the key to prove the blow-up result. We denote that our global existence and asymptotic behavior as $t\rightarrow\infty$  for small data solutions are obtained in the case when $1\leq n\leq 3,$ due to the nonlocal in time nonlinearity. While the blow-up result is done in any dimensional space. Let us present our main results. 

First, the following local well-posedness result is needed. 
\begin{proposition}\label{prop1}
Let $1<p\leq {n}/{(n-2)}$ for $n\geq 3,$ and $p\in(1,\infty)$ for $n=1,2.$ Under the assumptions $(\ref{condition1})$-$(\ref{condition2})$ and $\gamma\in(0,1),$ the problem $(\ref{eq1})$ possesses a unique maximal mild solution $u,$ i.e. satisfies the integral equation $(\ref{mildsolution})$ below, such that
$$u\in C([0,T_{\max}),H^1(\mathbb{R}^n))\cap C^1([0,T_{\max}), L^2(\mathbb{R}^n)),$$
where $0< T_{\max}\leq\infty.$ Moreover, $u(t,\cdotp)$ is supported in the ball $B(t+K).$ In addition:
\begin{equation}\label{alternative}
\mbox{either}\;\;T_{\max}=\infty \quad\mbox{or else}\quad T_{\max}<\infty \;\;\mbox{and}\;\; \|u(t)\|_{H^1}+\|u_t(t)\|_2\rightarrow\infty\;\;\mbox{as}\;\; t\rightarrow T_{\max}.
 \end{equation}
\end{proposition}
\begin{rmk}
We say that $u$ is a global solution of $(\ref{eq1})$ if $T_{\max}=\infty,$ while in the case of $T_{\max}<\infty,$ we say that $u$ blows up in finite time.
\end{rmk}
\noindent Now, set 
$$p_\gamma:=1+\frac{2(2-\gamma)}{(n-2+2\gamma)_+},\quad p_1:=1+\frac{2(3-2\gamma)}{(n-2+2\gamma)_+},\quad p_2:=1+\frac{4(3-2\gamma)}{(n-4+4\gamma)_+}\quad\hbox{and}\quad p_3:=1+\frac{n+2(5-4\gamma)}{(n-2+4\gamma)_+}.$$
As 
$$(p_\gamma={n}/{(n-2)={1}/{\gamma}})\Longleftrightarrow(\gamma={(n-2)}/{n}),$$ 
this imply, in the case when ${(n-2)}/{n}<\gamma,$ that $p_\gamma=\max\{{1}/{\gamma}\;;\;p_\gamma\}< {n}/{(n-2)}.$ Moreover, $p_\gamma<\min_{1\leq n\leq3}(p_n).$\\
We note that
$$p_\gamma,p_1\rightarrow 1+2/n=p_c(n),\quad p_2\rightarrow (2\gamma+1)/(2\gamma-1)>p_c(2)\quad\hbox{and}\quad p_3\rightarrow 2>p_c(3)\qquad\hbox{as}\;\gamma\rightarrow 1.$$

Our global existence result is the following
\begin{theorem}\label{theo1}
 Let $1\leq n\leq3,$ $p>1,$ $\gamma\in({1}/{2},1)$ for $n=1,2$ and $\gamma\in({11}/{16},1)$ for $n=3.$ Assume that the initial data satisfy $(\ref{condition1})$-$(\ref{condition2})$ such that $\|u_0\|_{H^1}+ \|u_1\|_{L^2}$
is sufficiently small. If $p_n<p$  then the problem $(\ref{eq1})$ admits a unique global mild solution 
$$u\in C([0,\infty),H^1(\mathbb{R}^n))\cap C^1([0,\infty),L^2(\mathbb{R}^n)).$$
\end{theorem}
\noindent Note that, the requirement $\gamma\in({11}/{16},1)$ is just to assure that $p_3<n/(n-2)$ when $n=3.$

The second result is the finite time blow-up of the solution under some positive data which shows that the assumption on the exponent in the above theorem (for $n=1$ and $\gamma\rightarrow1$) is critical and it is exactly the same critical exponent to the semilinear heat equation $u_t-\Delta u =|u|^p.$ Moreover, we conjecture that $p_1$ will be the critical exponent of $(\ref{eq1})$ which is the critical one to the corresponding semilinear heat equation $u_t-\Delta u =\int_0^t(t-s)^{-\gamma}|u(s)|^p\,ds$ founded  by Cazenave, Dickstein and Weissler \cite{CDW} and Fino and Kirane \cite{Finokirane}.\\

\begin{theorem}\label{theo2}${}$\\
\noindent $i)\;$ Let $1<p\leq {n}/{(n-2)}$ for $n\geq 3,$ and $p\in(1,\infty)$ for $n=1,2.$ Assume that ${(n-2)}/{n}<\gamma<1$ and $(u_0,u_1)$ satisfy $(\ref{condition1})$-$(\ref{condition2})$  such that
\begin{equation}\label{condition3}
\int_{\mathbb{R}^n}u_i(x)\,dx>0,\quad i=0,1.
\end{equation}
If $p\leq p_\gamma,$ then the mild solution of the problem $(\ref{eq1})$ blows up in finite time.\\
\noindent $ii)\;$ Let $n\geq 3$ and $1<p\leq {n}/{(n-2)}.$ Assume that $\gamma\leq{(n-2)}/{n}$ and $(u_0,u_1)$ satisfy $(\ref{condition1})$ and $(\ref{condition3}),$ then the mild solution of the problem $(\ref{eq1})$ blows up in finite time.
\end{theorem}
As the by-product of our analysis in Theorem \ref{theo1}, we have the following result concerning the asymptotic behavior as $t\rightarrow\infty$ of solutions.
\begin{theorem}\label{theo3}
Under the assumptions of Theorem $\ref{theo1},$ the asymptotic behavior of the small data global solution $u$ of $(\ref{eq1})$ is given by
\begin{equation}\label{asympbehavior} 
\|Du(t,\cdotp)\|_{L^2(\mathbb{R}^{n}\setminus B(t^{{1/2}+\delta}))}=\mathcal{O}(e^{-t^{2{\delta}/{4}}}),\quad t\rightarrow\infty,
\end{equation}
that is the solution decays exponentially outside every ball $B(t^{1/2+\delta}),\delta>0.$ Moreover, the total energy satisfies
\begin{equation}\label{decayenergie} 
\|Du(t,\cdotp)\|_{L^2(\mathbb{R}^n)}=\mathcal{O}(t^{-n/4+1/2-\gamma}),\quad t\rightarrow\infty,
\end{equation}
for $n=1,$
\begin{equation}\label{decayenergie1} 
\|Du(t,\cdotp)\|_{L^2(\mathbb{R}^n)}=\mathcal{O}(t^{1/2-\gamma}),\quad t\rightarrow\infty,
\end{equation}
for $n=2$ and
\begin{equation}\label{decayenergie2} 
\|Du(t,\cdotp)\|_{L^2(\mathbb{R}^n)}=\mathcal{O}(t^{-\gamma}),\quad t\rightarrow\infty,
\end{equation}
for $n=3.$
\end{theorem}
As we have seen, we are restricted ourselves in the case of compactly supported data. This restriction leads us to the finite propagation speed property of the wave which plays an important role in the proof of the global solvability. The blow-up result and the local existence theorem could be proved removing the requirement for the compactness assumptions on the support of the initial data. For the global existence without assuming the compactness of support on the initial data, we refer the reader to \cite{Hosonoogawa,Ikehata,Ikehatatanizawa,Narazaki,Nishihara} where we have to take $u_0\in H^1(\mathbb{R}^n)\cap L^1(\mathbb{R}^n)$ and $u_1\in L^2(\mathbb{R}^n)\cap L^1(\mathbb{R}^n).$
\begin{rmk}
It is still open to show corresponding global existence of solutions, with small initial data, for $p_\gamma<p<p_n$ $(1\leq n\leq 3)$ and for $p_\gamma<p$ $(n\geq4).$
\end{rmk}

This paper is organized as follows: in Section \ref{sec2}, we present some definitions and properties concerning the fractional integrals and derivatives. Section \ref{sec3} contains the proofs of the global existence theorem (Theorem \ref{theo1}) and the asymptotic behavior of solution (Theorem \ref{theo3}). Section \ref{sec4} is devoted to the proof of the blow-up result (Theorem \ref{theo2}). Finally, to make this paper self-contained, we shall sketch the proof of the local existence of solution (Proposition \ref{prop1}) in \ref{appendix}.


\section{Preliminaries}\label{sec2}
\setcounter{equation}{0}
In this section, we give some preliminary properties on the fractional integrals and fractional derivatives that will be used in the proof of Theorem \ref{theo2}.\\
If $AC[0,T]$ is the space of all functions which are
absolutely continuous on $[0,T]$ with $0<T<\infty,$ then, for $f\in
AC[0,T],$ the left-handed and right-handed Riemann-Liouville
fractional derivatives $D^{\alpha}_{0|t}f(t)$ and
$D^{\alpha}_{t|T}f(t)$ of order $\alpha\in(0,1)$ are defined by 
\begin{equation}\label{P1}
D^\alpha_{0|t}f(t):=\partial_tJ^{1-\alpha}_{0|t}f(t)\quad\mbox{and}\quad D^\alpha_{t|T}f(t):=-\frac{1}{\Gamma(1-\alpha)}\partial_t\int_t^T(s-t)^{-\alpha}f(s)\,ds,\quad t\in[0,T],\end{equation}
where 
\begin{equation}\label{P2}
J^\alpha_{0|t}g(t):=\frac{1}{\Gamma(\alpha)}\int_0^t(t-s)^{\alpha-1}g(s)\,ds
\end{equation}
is the Riemann-Liouville fractional integral, for all $g\in L^q(0,T)$ $(1\leq q\leq\infty).$ We refer the reader to \cite{Kilbas} for the definitions above. Furthermore, for every $f,g\in C([0,T])$ such that
$D^\alpha_{0|t}f(t),D^\alpha_{t|T}g(t)$ exist and are continuous,
for all $t\in[0,T],$ $0<\alpha<1,$ we have the formula of
integration by parts (see  $(2.64)$ p. 46 in \cite{SKM})
\begin{equation}\label{P3}
\int_0^T \left(D^\alpha_{0|t}f\right)(t)g(t)\,dt \;=\; \int_0^T
f(t)\left(D^\alpha_{t|T}g\right)(t)\,dt.
\end{equation}
Note also that, for all $f\in AC^{n+1}[0,T]$ and all integer $n\geq0,$ we have (see $(2.2.30)$ in
\cite{Kilbas})
\begin{equation}\label{P4}
(-1)^{n}\partial_t^n.D^\alpha_{t|T}f=D^{n+\alpha}_{t|T}f,
\end{equation}
where
$$AC^{n+1}[0,T]:=\left\{f:[0,T]\rightarrow\mathbb{R}\;\hbox{and}\;\partial_t^nf\in
AC[0,T]\right\}$$
and $\partial_t^n$ is the usual $n$ times derivative. Moreover, for all $1\leq q\leq\infty,$ the
following formula (see \cite[Lemma 2.4 p.74]{Kilbas})
\begin{equation}\label{P5}
    D^\alpha_{0|t}J^\alpha_{0|t}=Id_{L^q(0,T)}
\end{equation}
holds almost everywhere on $[0,T].$\\
In the proof of Theorem \ref{theo2}, the following results are useful: if $w_1(t)=\left(1-{t}/{T}\right)_+^\sigma,$
$t\geq0,$ $T>0,$ $\sigma\gg1,$ then
\begin{equation}\label{P6}
D_{t|T}^\alpha w_1(t)=CT^{-\sigma}(T-t)_+^{\sigma-\alpha},\quad D_{t|T}^{\alpha+1} w_1(t)=CT^{-\sigma}(T-t)_+^{\sigma-\alpha-1},\quad D_{t|T}^{\alpha+2} w_1(t)=CT^{-\sigma}(T-t)_+^{\sigma-\alpha-2},
\end{equation}
for all $\alpha\in(0,1);$ so
\begin{equation}\label{P7}
\left(D_{t|T}^\alpha w_1\right)(T)=0,\quad \left(D_{t|T}^\alpha w_1\right)(0)=C\;T^{-\alpha},\quad \left(D_{t|T}^{\alpha+1} w_1\right)(T)=0\quad\mbox{and}\quad\left(D_{t|T}^{\alpha+1} w_1\right)(0)=C\;T^{-\alpha-1}.
\end{equation}
For the proof of this results, see \cite[Preliminaries]{Finogeorgiev}. Furthermore, the following lemma is useful to prove Theorem \ref{theo1}.
\begin{lemma}\label{singularity}$(\mbox{\cite[Lemma 4.1]{Cui}})\;$
Suppose that $0\leq\theta<1,$ $a\geq0$ and $b\geq0.$ Then there exists a constant $C>0$ depending only on $a,b$ and $\theta$ such that fot all $t>0,$
$$
\int_0^t(t-\tau)^{-\theta}(1+t-\tau)^{-a}(1+\tau)^{-b}\,d\tau\leq \left\{\begin{array}{ll}
\,\,\displaystyle{C(1+t)^{-\min(a+\theta,b)}}&\displaystyle{\hbox{if}\;\max(a+\theta,b)>1,}\\
{}\\
\displaystyle{C(1+t)^{-\min(a+\theta,b)}\ln(2+t)}&\displaystyle{\hbox{if}\;\max(a+\theta,b)=1,}\\
{}\\
\displaystyle{C(1+t)^{1-a-\theta-b}}&\displaystyle{\hbox{if}\;\max(a+\theta,b)<1.}\\
\end{array}
\right.
$$
\end{lemma}
Throughout this paper, positive constants will be denoted by $C$ and will change from line to line. 


\section{Global existence and asymptotic behavior}\label{sec3}
In view of the Proposition \ref{prop1}, global existence of a solution follows from the boundedness of its energy at all times. To obtain such a priori estimates, we shall proceed our proof based on the weighted energy method recently developed in Todorova and Yordanov \cite{Todorovayordanov}. We begin by defining
\begin{equation}\label{weightfunction}
\psi(x,t)=\frac{1}{2}(t+K-\sqrt{(t+K)^2-|x|^2}),\quad |x|<t+K.
\end{equation}
It is easily checked that $\psi_t<0,$
\begin{equation}\label{propertiespsi}
0<\psi(x,t)<\frac{K}{2}
\end{equation}
and, since
$$\sqrt{(t+K)^2-|x|^2}\leq t+K-{|x|^2}/{[2(t+K)]},$$
the function $\psi$ satisfies the inequality
\begin{equation}\label{inequoverpsi}
\psi(x,t)\geq \frac{|x|^2}{4(t+K)}.
\end{equation} 

\noindent{\bf Proof of Theorem \ref{theo1}.} Let $u$ be the local solution of the problem $(\ref{eq1})$ in $[0,T_{\max}).$ Let us introduce the energy functional
\begin{equation}\label{functional}
W(t):=(1+t)^j\|Du(t,\cdotp)\|_2, 
\end{equation}
where 
$$j:=n/4-1/2+\gamma\quad (n=1),\qquad j:=\gamma-1/2\quad (n=2)\qquad\hbox{and}\qquad  j:=\gamma\quad (n=3).$$
We will show that $W(t)\leq CI_0,$ where $I_0:=\|u_0\|_{H^1}+\|u_1\|_2$ is small enough. This not only gives the global existence but also shows that, for $n=1$ and $\gamma\rightarrow1,$ the solution decays at least as fast as that of the linear part $u_{tt}-\Delta u+u_t=0.$ For the rate of the linear problem, see $(\ref{linearpartw})$ below.

The estimate $(\ref{functional})$ will be done by the following lemmas. 
\begin{lemma}\label{firstestimation}
Let $1\leq n\leq3,$ $\gamma\in(1/2,1)$ for $n=1,2$ and $\gamma\in(11/16,1)$ for $n=3.$ For all $\delta>0$ and all $t\in[0,T_{\max}),$ the following weighted energy estimate holds
\begin{equation}\label{newestimation1}
(1+t)^j\|Du(t,\cdotp)\|_2\leq CI_0+C(\max_{[0,t]}(1+\tau)^\beta\|e^{\delta\psi(\tau,\cdotp)}u(\tau,\cdotp)\|_{2p})^p,
\end{equation}
where $\beta>{n}/{4p}+{(2-\gamma})/{p}$ for $n=1,3$ and $\beta>{(2-\gamma})/{p}$ for $n=2.$
\end{lemma}
\begin{lemma}\label{newgagliardo}$(\mbox{\cite[Proposition 2.4]{Todorovayordanov}})$
Let $\theta(q)=n(1/2-1/q)$ and $0\leq \theta(q)\leq 1,$ and let $0<\sigma\leq1.$ If $u\in H^1(\mathbb{R}^n)$ with supp$u\subset B(t+K),$ $t\geq0.$ Then
\begin{equation}\label{neweqgagliardo}
\|e^{\sigma\psi(t,\cdotp)}u\|_q\leq C_K(1+t)^{(1-\theta(q))/2}\|\nabla u\|^{1-\sigma}_2\|e^{\psi(t,\cdotp)}\nabla u\|_2^\sigma,.
\end{equation}
where $\psi(t,x)$ is the weight function from $(\ref{weightfunction}).$
\end{lemma}

\noindent We postpone the proof of Lemma \ref{firstestimation} to the end of this section.\\
It follows from Lemma \ref{firstestimation} that
\begin{equation}\label{esti0}
W(t)\leq CI_0+C(\max_{[0,t]}(1+\tau)^\beta\|e^{\delta\psi(\tau,\cdotp)}u(\tau,\cdotp)\|_{2p})^p.
\end{equation}
On the other hand, Lemma \ref{newgagliardo} with $q=2p$ and $\sigma=\delta\leq1$ gives
\begin{eqnarray}\label{esti1}
\|e^{\delta\psi(\tau,\cdotp)}u(\tau,\cdotp)\|_{2p}&\leq& C(1+\tau)^{1-\theta(2p))/2}\|\nabla u\|^{1-\delta}_2\|e^{\psi(t,\cdotp)}\nabla u\|_2^\delta\nonumber\\
&\leq&C(1+\tau)^{(1-\theta(2p))/2-j}W(\tau),
\end{eqnarray}
where we have used (\ref{propertiespsi}).\\
Using $(\ref{esti1}),$ we obtain from $(\ref{esti0})$
\begin{equation}\label{esti3}
W(t)\leq CI_0+C\left(\max_{[0,t]}(1+\tau)^{\beta+(1-\theta(2p))/2-j}W(\tau)\right)^{p}.
\end{equation}
Set $\beta={n}/{4p}+{(2-\gamma)}/{p}+\nu$ for $n=1,3$ and $\beta={(2-\gamma)}/{p}+\nu$ for $n=2,$ $\nu>0,$ then if we compute the exponent of $(\tau+1)$ in the right side of $(\ref{esti3}),$ we obtain
\begin{equation}\label{quantity1}
\beta+(1-\theta(2p))/2-j=\left\{\begin{array}{ll}
\,\,\displaystyle{\nu-\frac{n}{2p}\left[p(1-{2(1-\gamma)}/{n})-1-{2(2-\gamma)}/{n}\right],}&\displaystyle{\hbox{if}\;n=1,}\\
{}\\
\displaystyle{\nu-\frac{n}{4p}\left[p(1-{4(1-\gamma)}/{n})-1-{4(2-\gamma)}/{n}\right],}&\displaystyle{\hbox{if}\;n=2,}\\
{}\\
\displaystyle{\nu-\frac{n}{4p}\left[p(1+{2(2\gamma-1)}/{n})-2-{4(2-\gamma)}/{n}\right],}&\displaystyle{\hbox{if}\;n=3.}\\
\end{array}
\right.
\end{equation}
As $p>p_n,$ we deduce, choosing $\nu$ small enough, that the quantities in $(\ref{quantity1})$ are negative. Hence, we can rewrite $(\ref{esti3})$ like
\begin{equation}\label{esti4}
\max_{[0,t]}W(\tau)\leq CI_0+C(\max_{[0,t]}W(\tau))^{p}.
\end{equation}
Now, write $I_0=\|u_0\|_{H^1}+\|u_1\|_2=C\varepsilon,$ for small $\varepsilon>0$ which is determined later, and put
$$T^*=\sup\{t\geq0:\;W(t)\leq 2C\varepsilon\}.$$
Then, $(\ref{esti4})$ implies $W(t)\leq C\varepsilon+C\varepsilon^{p}.$ Therefore, taking small $\varepsilon$ such that $C\varepsilon+C\varepsilon^{p}<2C\varepsilon$ we conclude that $T^*=\infty$ ( For details we refer the reader to \cite[Proposition 2.1]{Ikehatatanizawa} and \cite[Proposition 2.1]{Nishiharazhao}), i.e.
\begin{equation}\label{esti5}
W(t)=(1+t)^{j}\|Du(t,\cdotp)\|_2\leq C\varepsilon,\quad t\geq0.
\end{equation}
Thus we have completed the proof of Theorem \ref{theo1}.$\hfill\square$\\

\noindent{\bf Proof of Theorem \ref{theo3}.} The estimate $(\ref{decayenergie})-(\ref{decayenergie2})$ follows directly from $(\ref{esti5}).$ Next, it follows from inequality 
$(\ref{propertiespsi})$-$(\ref{inequoverpsi})$ and estimate $(\ref{esti5})$ that
$$C\varepsilon\geq \|e^{\psi(t,\cdotp)}Du(t,\cdotp)\|_{L^2(\mathbb{R}^n)}\geq \|e^{{|\cdotp|^2}/{4(t+K)}}Du(t,\cdotp)\|_{L^2(\mathbb{R}^n\setminus B(t^{1/2+\delta}))}\geq e^{{t^{1+2\delta}}/{4(t+K)}} \|Du(t,\cdotp)\|_{L^2(\mathbb{R}^n\setminus B(t^{1/2+\delta}))},$$
where we have used the fact that $j>0,$ which implies $(\ref{asympbehavior}).$ $\hfill\square$\\

To show Lemma $\ref{firstestimation},$ we need a linear estimates for the fundamental solution of the following linear damped wave equation
\begin{equation}\label{linearequation}
w_{tt}-\Delta w+w_t=0,\quad w(0,x)=u_0(x),\quad w_t(0,x)=u_1(x),
\end{equation}
for $t\in(0,\infty)\times\mathbb{R}^n.$ Let $K_0(t),K_1(t)$ be
$$K_0(t):=e^{-\frac{t}{2}}\cos\{ta(|\nabla|)\},\qquad K_1(t):=e^{-\frac{t}{2}}\frac{\sin\{ta(|\nabla|)\}}{a(|\nabla|)},$$
where
$$\mathcal{F}[a(|\nabla|)](\xi)=a(\xi)=\left\{\begin{array}{ll}
\,\,\displaystyle{\sqrt{|\xi|^2-1/4},}&\displaystyle{|\xi|>1/2,}\\
{}\\
\displaystyle{i\sqrt{1/4-|\xi|^2},}&\displaystyle{|\xi|\leq 1/2.}\\
\end{array}
\right.
$$
Note that $K_0(t)+{1}/{2}K_1(t)=\partial_t K_1(t).$ Then the solution of $(\ref{linearequation})$  is given (cf. \cite{Matsumura}) through the Fourier transform by $K_0(t)$ and $K_1(t)$ as
\begin{equation}\label{homogeneoussolution}
w(t,x)=K_0(t)\ast u_0+K_1(t)\ast\left(\frac{1}{2}u_0+u_1\right).
\end{equation}
The Duhamel principle implies that the solution $u(t,x)$ of nonlinear equation $(\ref{eq1})$ solves the integral equation
\begin{equation}\label{mildsolution}
u(t,x)=w(t,x)+\Gamma(\alpha)\int_0^tK_1(t-\tau)\ast J^{\alpha}_{0|\tau}(|u|^p)(\tau)\,d\tau,
\end{equation}
where $\alpha:=1-\gamma$ and $J^{\alpha}_{0|t}$ is given by $(\ref{P2}).$ We can now state Matsumura's result, on the estimate of $K_0(t)$ and $K_1(t),$ as follows:
\begin{lemma}\label{matsumuraresult}$(\mbox{\cite{Matsumura}})\;$ If $f\in L^m(\mathbb{R}^n)\cap H^{k+|\nu|-1}(\mathbb{R}^n)$ $(1\leq m\leq2),$ then
$$\|\partial_t^k\nabla_x^\nu K_1(t)\ast f\|_2\leq C(1+t)^{n/4-n/(2m)-|\nu|/2-k}(\|f\|_m+\|f\|_{H^{k+|\nu|-1}(\mathbb{R}^n)}).$$
\end{lemma}

\noindent{\bf Proof of Lemma \ref{firstestimation}.} We begin to estimate the linear term $\|Dw(t,\cdotp)\|_2.$ It is not difficult to see, using Lemma \ref{matsumuraresult} with $m=1,$ that
\begin{equation}\label{linearpartw}
\|Dw(t,\cdotp)\|_2\leq C(1+t)^{-n/4-1/2}(\|u_0\|_{H^1}+\|u_0\|_1+\|u_1\|_2+\|u_1\|_1)\leq CI_0(1+t)^{-n/4-1/2}\leq  CI_0(1+t)^{-j}.
\end{equation}
To estimate the nonlinear term in $(\ref{mildsolution}),$ we have to distinguish two cases:\\

\noindent $\bullet$ \underline{Case of $n=1,3$}: Apply Lemma $\ref{matsumuraresult}$ with $m=1$ to get 
\begin{eqnarray}\label{esti6}
I:=\int_0^t\| DK_1(t-\tau)\ast J^{\alpha}_{0|\tau}(|u|^p)(\tau)\|_2\,d\tau&\leq&C\int_0^t(t-\tau+1)^{-n/4-1/2}\left(\|J^{\alpha}_{0|\tau}(|u|^p)(\tau)\|_1+\|J^{\alpha}_{0|\tau}(|u|^p)(\tau)\|_2\right)\,d\tau\nonumber\\
&\leq&C\int_0^t(t-\tau+1)^{-n/4-1/2}\left(J^{\alpha}_{0|\tau}\|u(\tau)\|_p^p+J^{\alpha}_{0|\tau}\|u(\tau)\|_{2p}^p\right)\,d\tau.
\end{eqnarray}
To transform the $L^p$-norm into a weighted $L^{2p}$-norm, we use the Cauchy inequality
\begin{eqnarray*}
\|u(\tau,\cdotp)\|_p^p&\equiv&\int_{B(\tau+K)}|u(\tau,x)|^p\,dx\\
&\leq&\left(\int_{B(\tau+K)}e^{-2p\delta\psi(\tau,x)}\,dx\right)^{1/2}\left(\int_{B(\tau+K)}e^{2p\delta\psi(\tau,x)}|u(\tau,x)|^{2p}\,dx\right)^{1/2},
\end{eqnarray*}
for $\delta>0.$ From $(\ref{inequoverpsi}),$ we have $\psi(\tau,x)\geq{|x|^2}/{4(\tau+K)}$ for $x\in B(\tau+K),$ so the first integral is estimated as follows
$$\int_{B(\tau+K)}e^{-2p\delta\psi(\tau,x)}\,dx\leq \int_{B(\tau+K)}e^{-p\delta{|x|^2}/{2(\tau+k)}}\,dx\leq  \int_{\mathbb{R}^n}e^{-p\delta{|x|^2}/{2(\tau+k)}}\,dx\equiv \left(\frac{2\pi}{p\delta}\right)^{n/2}(\tau+K)^{n/2}.
$$
Thus, for the norm $\|u(\tau,\cdotp)\|_p$ in $(\ref{esti6})$ we obtain the weighted estimate
\begin{equation}
\|u(\tau,\cdotp)\|_p^p\leq C_{K,\delta}(\tau+1)^{n/4}\|e^{\delta\psi(\tau,\cdotp)}u(\tau,\cdotp)\|_{2p}^p,\quad \delta>0.
\end{equation}\label{esti7}
Next, as $\psi>0,$  the norm $\|u(\tau,\cdotp)\|_{2p}$ in $(\ref{esti6})$ can obviously be estimated by
\begin{equation}\label{esti8}
\|u(\tau,\cdotp)\|_{2p}^p\leq C_\delta(\tau+1)^{n/4}\|e^{\delta\psi(\tau,\cdotp)}u(\tau,\cdotp)\|_{2p}^p.
\end{equation}
Combining $(\ref{esti6})-(\ref{esti8}),$ we obtain
\begin{eqnarray*}\label{nonlinearpart1}
I&\leq &C \int_0^{t}(t-\tau+1)^{-n/4-1/2}\int_0^\tau(\tau-\sigma)^{-\gamma}((\sigma+1)^{n/(4p)}\|e^{\delta\psi(\sigma,\cdotp)}u(\sigma,\cdotp)\|_{2p})^p\,d\sigma\,d\tau\nonumber\\
&\leq&C(\max_{[0,t]}(\tau+1)^{\beta}\|e^{\delta\psi(\tau,\cdotp)}u(\tau,\cdotp)\|_{2p})^p \int_0^{t}(t-\tau+1)^{-n/4-1/2}\int_0^\tau(\tau-\sigma)^{-\gamma}(1+\sigma)^{-2(2-\gamma)}\,d\sigma\,d\tau.
\end{eqnarray*}
Using Lemma \ref{singularity}, we conclude that
\begin{equation}\label{nonlinearpart4}
I\leq C(1+t)^{-j} (\max_{[0,t]}(\tau+1)^{\beta}\|e^{\delta\psi(\tau,\cdotp)}u(\tau,\cdotp)\|_{2p})^p.
\end{equation}
Combining $(\ref{linearpartw})$ and $(\ref{nonlinearpart4}),$ we obtain $(\ref{newestimation1}).$ This complete the proof for $n=1,3.$\\

\noindent $\bullet$ \underline{Case of $n=2$}:  Apply here Lemma $\ref{matsumuraresult}$ with $m=2,$ we obtain
\begin{eqnarray*}\label{esti06}
J:=\int_0^t\| DK_1(t-\tau)\ast J^{\alpha}_{0|\tau}(|u|^p)(\tau)\|_2\,d\tau&\leq& C\int_0^t(t-\tau+1)^{-1/2}\|J^{\alpha}_{0|\tau}(|u|^p)(\tau)\|_2\,d\tau\nonumber\\
&\leq& C\int_0^t(t-\tau+1)^{-1/2}\int_0^\tau(\tau-\sigma)^{-\gamma}\|u(\sigma)\|_{2p}^p\,d\sigma\,d\tau.
\end{eqnarray*}
Then
\begin{eqnarray}\label{esti07}
J\leq C(\max_{[0,t]}(\tau+1)^{\beta}\|e^{\delta\psi(\tau,\cdotp)}u(\tau,\cdotp)\|_{2p})^p \int_0^{t}(t-\tau+1)^{-1/2}\int_0^\tau(\tau-\sigma)^{-\gamma}(1+\sigma)^{-2(2-\gamma)}\,d\sigma\,d\tau.
\end{eqnarray}
By Lemma \ref{singularity}, $(\ref{esti07})$ implies
\begin{equation}\label{nonlinearpart40}
J\leq C(1+t)^{-j} (\max_{[0,t]}(\tau+1)^{\beta}\|e^{\delta\psi(\tau,\cdotp)}u(\tau,\cdotp)\|_{2p})^p.
\end{equation}
Combining $(\ref{linearpartw})$ and $(\ref{nonlinearpart40}),$ we obtain $(\ref{newestimation1}).$ This complete the proof for $n=2.$ $\hfill\square$\\


\section{Blow-up result}\label{sec4}
\setcounter{equation}{0} 
In this section we devote ourselves to the proof of Theorem \ref{theo2}.  We start by introducing the definition of the weak solution of $(\ref{eq1}).$
\begin{definition}$(\mbox{Weak solution})\;$ Let $T>0,$ $\gamma\in(0,1)$ and $u_0,u_1\in L_{loc}^1(\mathbb{R}^n).$ We say that $u$ is a weak solution if $u\in L^p((0,T),L_{loc}^p(\mathbb{R}^n))$ and satisfies
\begin{eqnarray}\label{weaksolution}
&{}&\Gamma(\alpha)\int_0^T\int_{\mathbb{R}^n}J^{\alpha}_{0|t}(|u|^{p})\varphi\,dx\,dt+\int_{\mathbb{R}^n}u_1(x)\varphi(0,x)\,dx+\int_{\mathbb{R}^n}u_0(x)(\varphi(0,x)-\varphi_t(0,x))\,dx\nonumber\\
&{}&=\int_0^T\int_{\mathbb{R}^n}u\varphi_{tt}\,dx\,dt-\int_0^T\int_{\mathbb{R}^n}u\varphi_t\,dx\,dt-\int_0^T\int_{\mathbb{R}^n}u\Delta\varphi\,dx\,dt,
\end{eqnarray}
for all compactly supported function $\varphi\in C^2([0,T]\times\mathbb{R}^n)$ such that $\varphi(\cdotp,T)=0$ and $\varphi_t(\cdotp,T)=0,$ where $\alpha=1-\gamma.$ 
\end{definition}
Next, the following lemma is useful for the proof of Theorem \ref{theo2}. The proof of this lemma is much the same procedure as in the proof of \cite[Lemma 2]{Finogeorgiev}.
\begin{lemma}\label{mildweak}$(\mbox{Mild $\rightarrow$ Weak})\;$ Let $T>0$ and $\gamma\in(0,1).$ Suppose that $1<p\leq {n}/{(n-2)},$ if $n\geq 3,$ and $p\in(1,\infty),$ if $n=1,2.$ If $u\in C([0,T],H^1(\mathbb{R}^n))\cap C^1([0,T],L^2(\mathbb{R}^n))$ is the mild solution of $(\ref{eq1}),$ then $u$ is a weak solution of $(\ref{eq1}).$
\end{lemma}
 \noindent{\bf Remark.}$\;$ We need the mild solution to use, in the proof of Theorem \ref{theo2}, the alternative $(\ref{alternative}).$ Without this properties, we say that we have a nonexistence of global solution and not a blow-up result.\\
 
 \noindent{\bf Proof of Theorem \ref{theo2}.}$\;$ We assume on the contrary, using $(\ref{alternative}),$ that $u$ is a global mild solution of $(\ref{eq1}).$ So, from Lemma \ref{mildweak} we have
 \begin{eqnarray}\label{newweaksolution}
&{}&\Gamma(\alpha)\int_0^T\int_{\mbox{supp$\varphi$}}J^{\alpha}_{0|t}(|u|^{p})\;\varphi\,dx\,dt+\int_{\mbox{supp$\varphi$}}u_1(x)\varphi(0,x)\,dx+\int_{\mbox{supp$\varphi$}}u_0(x)(\varphi(0,x)-\varphi_t(0,x))\,dx\nonumber\\
&{}&=\int_0^T\int_{\mbox{supp$\varphi$}}u\;\varphi_{tt}\,dx\,dt-\int_0^T\int_{\mbox{supp$\varphi$}}u\;\varphi_t\,dx\,dt-\int_0^T\int_{\mbox{supp$\Delta\varphi$}}u\;\Delta\varphi\,dx\,dt,
\end{eqnarray}
for all $T>0$ and all compactly supported test function $\varphi\in C^2([0,T]\times\mathbb{R}^n)$ such that $\varphi(\cdotp,T)=0$ and $\varphi_t(\cdotp,T)=0,$ where $\alpha=1-\gamma.$ Let $\varphi(x,t)=D^\alpha_{t|T}\left(\tilde{\varphi}(x,t)\right):=D^\alpha_{t|T}\left(\varphi^\ell_1(x)\varphi_2(t)\right)$ with
$\varphi_1(x):=\Phi\left({|x|}/{B}\right),$
$\varphi_2(t):=\left(1-{t}/{T}\right)^\eta_+,$ where $D^\alpha_{t|T}$ is given by $(\ref{P1}),$ $\ell,\eta\gg1$ and $\Phi\in C^\infty(\mathbb{R}_+)$ be a cut-off non-increasing function such that
$$\Phi(r)=\left\{\begin {array}{ll}\displaystyle{1}&\displaystyle{\quad\mbox{if }0\leq r\leq 1}\\
\displaystyle{0}&\displaystyle{\quad\mbox {if }r\geq 2,}
\end {array}\right.$$\\
$0\leq \Phi \leq 1$ and $|\Phi^{'}(r)|\leq C_1/r$ for all $r>0.$ The constant $B>0$ in the definition of $\varphi_1$ is fixed and will be chosen later. In the following, we denote by $\Omega(B)$ the support of $\varphi_1$ and by $\Delta(B)$ the set containing the support of  $\Delta\varphi_1$ which are defined as follows:
$$\Omega(B)=\{x\in\mathbb{R}^n:\;|x|\leq 2B\},\quad \Delta(B)=\{x\in\mathbb{R}^n:\;B\leq|x|\leq 2B\}.$$  
We return to $(\ref{newweaksolution}),$ which actually reads
  \begin{eqnarray}\label{newweaksolution1}
&{}&\Gamma(\alpha)\int_0^T\int_{\Omega(B)}J^{\alpha}_{0|t}(|u|^{p})D^\nu_{t|T}\tilde{\varphi}\,dx\,dt+\int_{\Omega(B)}u_1(x)D^\alpha_{t|T}\tilde{\varphi}(0,x)\,dx+\int_{\Omega(B)}u_0(x)(D^\alpha_{t|T}\tilde{\varphi}(0,x)-\partial_tD^\alpha_{t|T}\tilde{\varphi}(0,x))\,dx\nonumber\\
&{}&=\int_0^T\int_{\Omega(B)}u\;\partial_t^2D^\alpha_{t|T}\tilde{\varphi}\,dx\,dt-\int_0^T\int_{\Omega(B)}u\;\partial_tD^\alpha_{t|T}\tilde{\varphi}\,dx\,dt-\int_0^T\int_{\Delta(B)}u\;\Delta D^\alpha_{t|T}\tilde{\varphi}\,dx\,dt.
\end{eqnarray} 
 From $(\ref{P3}),(\ref{P4})$ and $(\ref{P7}),$ we conclude that
  \begin{eqnarray}\label{newweaksolution2}
&{}&\int_0^T\int_{\Omega(B)}D^\alpha_{0|t}J^{\alpha}_{0|t}(|u|^{p})\tilde{\varphi}\,dx\,dt+C\;T^{-\alpha}\int_{\Omega(B)}u_1(x)\varphi^\ell_1(x)\,dx+C(T^{-\alpha}+T^{-\alpha-1})\int_{\Omega(B)}u_0(x)\varphi^\ell_1(x)\,dx\nonumber\\
&{}&=C\int_0^T\int_{\Omega(B)}u(D^{2+\alpha}_{t|T}\tilde{\varphi}+D^{1+\alpha}_{t|T}\tilde{\varphi})\,dx\,dt-C\int_0^T\int_{\Delta(B)}u\;\Delta(\varphi^\ell_1)D^\alpha_{t|T}\varphi_2\,dx\,dt,
\end{eqnarray}  
where $D^\alpha_{0|t}$ is defined in $(\ref{P1}).$ Moreover, using $(\ref{P5})$ and the fact that $(\ref{condition3})$ implies
$\int_{\Omega(B)}\varphi^\ell_1(x) u_i(x)\geq0, i=0,1,$ it follows
  \begin{eqnarray}\label{newweaksolution3}
\int_0^T\int_{\Omega(B)}|u|^{p}\tilde{\varphi}\,dx\,dt&\leq&C\int_0^T\int_{\Omega(B)}|u|\varphi_1^\ell (D^{2+\alpha}_{t|T}\varphi_2+D^{1+\alpha}_{t|T}\varphi_2)\,dx\,dt\nonumber\\
&{}&+\;C\int_0^T\int_{\Delta(B)}|u|\varphi_1^{\ell-2}(|\Delta\varphi_1|+|\nabla\varphi_1|^2)D^\alpha_{t|T}\varphi_2\,dx\,dt\nonumber\\
&=:&I_1+I_2,
\end{eqnarray}  
where we have used the formula $\Delta(\varphi_1^\ell)=\ell\varphi_1^{\ell-1}\Delta\varphi_1+\ell(\ell-1)\varphi_1^{\ell-2}|\nabla\varphi_1|^2$ and $\varphi_1\leq1.$ Next we observe that by introducing the term $\tilde{\varphi}^{1/p}\tilde{\varphi}^{-1/p}$ in the right side of $(\ref{newweaksolution3})$ and applying Young's inequality 
we have
\begin{equation}\label{conditionI1}
I_1\leq\frac{1}{2p}\int_0^T\int_{\Omega(B)}|u|^{p}\tilde{\varphi}\,dx\,dt+C\int_0^T\int_{\Omega(B)}\varphi_1^\ell\varphi_2^{-1/(p-1)} ((D^{2+\alpha}_{t|T}\varphi_2)^{p'}+(D^{1+\alpha}_{t|T}\varphi_2)^{p'})\,dx\,dt,
\end{equation}
where $p'={p}/{(p-1)}.$ Similarly,
\begin{equation}\label{conditionI2}
I_2\leq  \frac{1}{2p}\int_0^T\int_{\Omega(B)}|u|^{p}\tilde{\varphi}\,dx\,dt+C\int_0^T\int_{\Omega(B)}\varphi_1^{\ell-2p'}\varphi_2^{-1/(p-1)} (|\Delta\varphi_1|^{p'}+|\nabla\varphi_1|^{2p'})(D^\alpha_{t|T}\varphi_2)^{p'}\,dx\,dt.
\end{equation}
Combining $(\ref{conditionI1})$ and $(\ref{conditionI2}),$ it follows from $(\ref{newweaksolution3})$ that
 \begin{eqnarray}\label{newweaksolution4}
\int_0^T\int_{\Omega(B)}|u|^{p}\tilde{\varphi}\,dx\,dt&\leq&C\int_0^T\int_{\Omega(B)}\varphi_1^\ell\varphi_2^{-1/(p-1)} ((D^{2+\alpha}_{t|T}\varphi_2)^{p'}+(D^{1+\alpha}_{t|T}\varphi_2)^{p'})\,dx\,dt\nonumber\\
&{}&+\;C\int_0^T\int_{\Omega(B)}\varphi_1^{\ell-2p'}\varphi_2^{-1/(p-1)} (|\Delta\varphi_1|^{p'}+|\nabla\varphi_1|^{2p'})(D^\alpha_{t|T}\varphi_2)^{p'}\,dx\,dt.
\end{eqnarray}  
At this stage, to prove $i),$ we have to distinguishes 2 cases.\\ 

\noindent $\bullet$ \underline{Case of $p<p_{\gamma}$}: in this case, we take $B=T^{1/2}.$ So, using $(\ref{P6})$ and the change of variables: $s=T^{-1}t,$ $y=T^{-1/2}x,$ we get from $(\ref{newweaksolution4})$ that 
\begin{equation}\label{newweaksolution5}
\int_0^T\int_{\Omega(T^{1/2})}|u|^{p}\tilde{\varphi}\,dx\,dt\leq C(T^{-(\alpha+2)p'+n/2+1}+T^{-(\alpha+1)p'+n/2+1}),
\end{equation}  
where $C$ is independent of $T.$ Letting $T\rightarrow\infty$ in $(\ref{newweaksolution5}),$ thanks to $p<p_{\gamma}$ and the Lebesgue dominated convergence theorem, it is yielded that
$$\int_0^\infty\int_{\mathbb{R}^n}|u|^p\,dx\,dt=0,$$
which implies $u(x,t)=0$ for all $t$ and a.e. $x.$ This contradicts our assumption $(\ref{condition3}).$\\

\noindent $\bullet$ \underline{Case of $p=p_{\gamma}$}: let $B=R^{-1/2}T^{1/2},$ where $1\ll R<T$ is such that when $T\rightarrow\infty$ we don't have $R\rightarrow\infty$ at the same time. Moreover, from the last case and the fact that $p=p_{\gamma},$ there exist a positive constant $D$ independent of $T$ such that
$$\int_0^\infty\int_{\mathbb{R}^n}|u|^p\,dx\,dt\leq D,$$
which implies that
\begin{equation}\label{conditioninfini}
\int_0^T\int_{\Delta(R^{-1/2}T^{1/2})}|u|^p\tilde{\varphi}\,dx\,dt\rightarrow0\quad\mbox{as}\;\;T\rightarrow\infty.
\end{equation}
On the other hand, using H\"older's inequality instead of Young's one, we estimate the integral $I_2$ in $(\ref{newweaksolution3})$ as follows:
 \begin{equation}\label{holderconditionI2}
 I_2\leq C\left(\int_0^T\int_{\Delta(R^{-1/2}T^{1/2})}|u|^p\tilde{\varphi}\right)^{1/p}\left(\int_0^T\int_{\Omega(R^{-1/2}T^{1/2})}\varphi_1^{\ell-2p'}\varphi_2^{-{1}/{(p-1)}}(|\Delta\varphi_1|^{p'}+|\nabla\varphi_1|^{2p'})(D^\alpha_{t|T}\varphi_2)^{p'}\,dx\,dt\right)^{{1}/{p'}}.
  \end{equation}
Similarly to the last case, substituting $(\ref{conditionI1})$ and $(\ref{holderconditionI2})$ into $(\ref{newweaksolution3}),$ taking account of $p=p_\gamma$ and the scaled variable $s=T^{-1}t, \;y=R^{1/2}T^{-1/2}x,$ we get
 $$\int_0^T\int_{\Omega(R^{-1/2}T^{1/2})}|u|^p\,dx\,dt \leq C(T^{-p'}R^{-n/2}+R^{-n/2})+CR^{1-{n}/{(2p')}}\left(\int_0^T\int_{\Delta(R^{-1/2}T^{1/2})}|u|^p\tilde{\varphi}\right)^{1/p}.$$
 Letting $T\rightarrow\infty,$ using $(\ref{conditioninfini}),$ we get
  $$\int_0^\infty\int_{\mathbb{R}^N}|u|^p\,dx\,dt \leq CR^{-n/2},$$
  which implies a contradiction, when $R\rightarrow\infty,$ with $(\ref{condition3}).$ This completes the proof of Theorem $\ref{theo2},$ $i).$\\
  For the proof of $ii),$ we have two possibility.\\
  
\noindent $\bullet$ \underline{If $\gamma<{(n-2)}/{n}$}: let $B=R$ with the same $R$ introduced in the case $p=p_{\gamma}.$ Then, taking the scaled variables $s=T^{-1}t,\;y=R^{-1}x,$ it follows from $(\ref{newweaksolution4})$ that
$$\int_0^T\int_{\Omega(R)}|u|^{p}\tilde{\varphi}\,dx\,dt\leq CR^n(T^{-(2+\alpha)p'+1}+T^{-(1+\alpha)p'+1})+CR^{n-2p'}T^{-\alpha p'+1}.$$
As $\gamma<{(n-2)}/{n}$ implies $p\leq{n}/{n-2}<{1}/{\gamma},$ we get a contradiction with $(\ref{condition3})$ by letting the following limits: first $T\rightarrow\infty,$ next $R\rightarrow\infty.$\\

\noindent $\bullet$ \underline{If $\gamma={(n-2)}/{n}$}: we have $p\leq {n}/{(n-2)}={1}/{\gamma}=p_\gamma.$ Using the first two cases, we get the contradiction. This completes the proof of Theorem $\ref{theo2},$ $ii).$ $\hfill\square$\\


\appendix
\section{}\label{appendix}
In this appendix let us sketch the proof of Proposition $\ref{prop1}.$ Let us define a semigroup $S(t):H^1(\mathbb{R}^n)\times L^2(\mathbb{R}^n)\rightarrow H^1(\mathbb{R}^n)\times L^2(\mathbb{R}^n)$ by
$$S(t):\left[\begin{array}{l}
u_0\\
u_1\\
\end{array}
\right]\mapsto\left[\begin{array}{l}
w\\
w_t\\
\end{array}
\right],$$
where $w\in C([0,\infty),H^1(\mathbb{R}^n))\cap C^1([0,\infty),L^2(\mathbb{R}^n))$ is the linear solution of $(\ref{linearequation})$ given by $(\ref{homogeneoussolution}).$ So, view of $(\ref{mildsolution}),$ a mild solution of the nonlinear problem $(\ref{eq1})$ is equivalent to  following integral equation:
\begin{equation}\label{newmildsolution}
U(t)=S(t)U_0+\int_0^t S(t-s)F(s)\,ds,
\end{equation}
where 
$$U(t)=\left[\begin{array}{l}
u(t,\cdotp)\\
u_t(t,\cdotp)\\
\end{array}
\right],\quad U_0=\left[\begin{array}{l}
u_0\\
u_1\\
\end{array}
\right],\quad F(s)=\left[\begin{array}{l}
0\\
J^\alpha_{0|s}(|u|^p)(s)\\
\end{array}
\right].$$
It sufficient now to prove the local existence of a solution of $(\ref{newmildsolution})$ in $H^1(\mathbb{R}^n)\times L^2(\mathbb{R}^n).$ Let $T>0$ and consider the following Banach space
$$E:=\{U={}^t(u,\varv):\;(u,\varv)\in C([0,T],H^1(\mathbb{R}^n)\times L^2(\mathbb{R}^n)), \;\mbox{supp}u(t,\cdotp)\subset B(K+t)\;\mbox{and}\;\|U\|_E\leq CM\},$$
where 
$$\|U\|_E:=\|u\|_{C([0,T];H^1(\mathbb{R}^n))}+\|\varv\|_{C([0,T];L^2(\mathbb{R}^n))}\quad\mbox{and}\quad M:=\|u_0\|_{H^1}+\|u_1\|_2.$$
In order to use the Banach fixed point theorem, we introduce the following map $\Phi$ on $E$ defined by
$$\Phi[U](t):=S(t)U_0+\int_0^t S(t-s)F(s)\,ds.$$
Now, for $U=(u,\varv)\in E,$ we have 
$$\|J^{\alpha}_{0|t}(|u|^p)(t)\|_{2}\leq Ct^{1-\gamma}\|u(t,\cdotp)\|^p_{2p}\leq Ct^{1-\gamma}\|u(t,\cdotp)\|^p_{H^1}\leq Ct^{1-\gamma}\|U\|_E^p,\quad t\in[0,T],$$
where we have used the Sobolev imbedding $H^1(\mathbb{R}^n)\subset L^{2p}(\mathbb{R}^n).$ Next, using Matsumura's result (Lemma \ref{matsumuraresult}) with $m=2$ and the finite propagation speed phenomena, we deduce via the Banach fixed point theorem that there exists a local solution $U\in E$ on a small interval $[0,T]$ satisfies $(\ref{newmildsolution}).$ For details, we refer the reader to \cite[Theorem 3.2]{Finokirane} and \cite[Theorem 6]{Finogeorgiev}. By consequence, there exist a local solution $u\in C([0,T],H^1(\mathbb{R}^n))\cap C^1([0,T],L^2(\mathbb{R}^n))$ satisfies $(\ref{mildsolution})$ and supp$u(t,\cdotp)\subset B(t+K).$  However, since our equation $(\ref{eq1})$ is nonautonomous, we prefer apply Gronwall's inequality to get the uniqueness (cf. \cite[Theorem~3.1]{CDW}). Indeed, if $u,\varv\in C([0,T],H^1(\mathbb{R}^n))\cap C^1([0,T],L^2(\mathbb{R}^n))$ are two mild solutions (i.e. satisfy $(\ref{mildsolution})$) for some $T>0,$ we have
\begin{eqnarray}\label{appendix1}
\|u(t)-\varv(t)\|_{H^1}&\leq&C\int_0^t\|K_1(t-\tau)\ast J^{\alpha}_{0|\tau}(|u|^p-|\varv|^p)(\tau)\|_{H^1}\,d\tau\nonumber\\
&\leq&C\int_0^t(1+t-\tau)^{-1/2}\| J^{\alpha}_{0|\tau}(|u|^p-|\varv|^p)(\tau)\|_{2}\,d\tau\nonumber\\
&\leq&C\int_0^t\| J^{\alpha}_{0|\tau}(|u|^p-|\varv|^p)(\tau)\|_{2}\,d\tau,
\end{eqnarray}
where we have used again Matsumura's result (Lemma \ref{matsumuraresult}) with $m=2.$ As $||u|^p-|\varv|^p|\leq C|u-\varv|(|u|^p+|\varv|^p),$ so by H\"older's inequality $(\|ab\|_{2}\leq \|a\|_{2p}\|b\|_{2p'})$ with $p'={p}/{(p-1)}$ and Sobolev's imbedding $(H^1\subset L^{2p}),$ we obtain
\begin{eqnarray}\label{appendix2}
\int_0^t\| J^{\alpha}_{0|\tau}(|u|^p-|\varv|^p)(\tau)\|_{2}\,d\tau&\leq& C\int_0^t J^{\alpha}_{0|\tau}(\|u-\varv\|_{H^1}(\|u\|_{H^1}^{p-1}+\|\varv\|^{p-1}_{H^1}))(\tau)\,d\tau\nonumber\\
&\leq& C\int_0^t \int_0^\tau (\tau-s)^{-\gamma}\|u(s,\cdotp)-\varv(s,\cdotp)\|_{H^1}\,ds\,d\tau\nonumber\\
&=& C\int_0^t \int_s^t (\tau-s)^{-\gamma}\|u(s,\cdotp)-\varv(s,\cdotp)\|_{H^1}\,d\tau\,ds\nonumber\\
&=& C\int_0^t (t-s)^{1-\gamma}\|u(s,\cdotp)-\varv(s,\cdotp)\|_{H^1}\,ds.
\end{eqnarray}
Combining $(\ref{appendix1})$ and $(\ref{appendix2}),$ we get
$$\|u(t)-\varv(t)\|_{H^1}\leq C\int_0^t (t-s)^{1-\gamma}\|u(s,\cdotp)-\varv(s,\cdotp)\|_{H^1}\,ds.$$
Using Gronwall's inequality, it follows that $u(t)\equiv \varv(t).$ As a consequence of this uniqueness result, we can extend our solution $u$ on a maximal interval $[0,T_{\max}).$ Moreover, if $T_{\max}<\infty,$ then $\|u(t,\cdotp)\|_{H^1}+\|u_t(t,\cdotp)\|_2\rightarrow\infty$ as $t\rightarrow T_{\max}.$ For details, see \cite[Theorem 3.1]{CDW} and \cite[Theorem 3.2]{Finokirane}.




\section*{Acknowledgements}
The author was supported by the Lebanese Association for Scientific Research (LASeR). The author would like to express sincere gratitude to Professor Mokhtar Kirane for valuable discussion. He thanks Professors Vladimir Georgiev, Thierry Cazenave and Flavio Dickstein for their helful remarks. 

\bibliographystyle{elsarticle-num}



\end{document}